\documentclass{amsart}
\usepackage{amsmath,amssymb}
\newtheorem{theorem}{Theorem}[section]
\newtheorem{lemma}[theorem]{Lemma}
\newtheorem{proposition}[theorem]{Proposition}
\newtheorem{corollary}[theorem]{Corollary}

\newtheorem{remark}[theorem]{Remark}
\newtheorem{example}[theorem]{Example}

\begin{document}

\title{RD-flatness and RD-injectivity}
\author{F. Couchot}
\address{Universit\'e de Caen Basse-Normandie, CNRS UMR
  6139 LMNO,
F-14032 Caen cedex 5, France}
\email{francois.couchot@unicaen.fr}

\begin{abstract}
It is proved that every commutative ring whose RD-injective modules are
 $\Sigma$-RD-injective is the product of a pure semi-simple ring and  a finite ring. A complete characterization 
 of commutative rings for which each artinian (respectively simple) module is RD-injective, is given. These results
  can be obtained by using the properties of  RD-flat modules and RD-coflat modules which are respectively 
  the RD-relativization of flat modules and fp-injective modules. It is also shown that a commutative ring is perfect if and only if each RD-flat module is RD-projective.
\end{abstract}

\maketitle

RD-purity is an important example of relative purity. It is the first notion of purity (for submodules) that appeared in the mathematical literature. Another reason of this importance is that for some classes of rings, not necessarily commutative, RD-purity coincides with purity. In
\cite{PPR99} G. Puninski, M. Prest and P. Rothmaler studied these rings that they called RD-rings. In particular they proved that the classe of commutative RD-rings is exactly the classe of arithmetic rings, i.e the rings with a distributive lattice of ideals. The first section of this paper is devoted to the study of RD-flat modules and RD-coflat modules which are respectively the RD-relativization of flat modules and fp-injective modules. It is interesting to note that the RD-flat modules form an example of additive accessible category by Proposition~\ref{P:flat}.5. Such categories were studied by Crawley-Boevey in \cite{CrBo94}, using the term of locally finitely presented additive categories. This allows us to prove that a commutative ring is perfect if and only if each RD-flat module is RD-projective (Corollary~\ref{C:perf}). On the other hand, it is shown (Theorem~\ref{T:flat}) that a ring is RD if and only if any right module is RD-flat (or RD-coflat). So, we can say that each RD-ring is absolutely RD-flat.
 
Recall that a ring is said to be pure semi-simple if each right module is pure-injective. It is known that a ring $R$ is pure semi-simple if and only if each right pure-injective $R$-module is $\Sigma$-pure-injective. In section~\ref{S:hull} (Theorem~\ref{T:cofinj}) we proved the RD-variation of this result: any right RD-coflat $R$-module is RD-injective if and only if any right RD-injective $R$-module is $\Sigma$-RD-injective. Moreover, any commutative ring which enjoys these properties, is the product of a pure semi-simple ring and a finite ring.

By \cite[Observation~3(4) and Theorem~6]{HZ00} every
artinian module over a commutative ring is
$\Sigma$-pure-injective. However, 
\cite[Example 4.6]{Couc03} is an example of simple module
over a noetherian domain that fails to be RD-injective. If $R$ is a commutative artinian ring, it is proved   (Proposition~\ref{P:artring}) that $R$ is an
RD-injective module if and only if $R$ is quasi-Frobenius.
In 
section~\ref{S:inje} we give a characterization of commutative rings for which any simple module is RD-injective (Theorem~\ref{T:glo}) and a characterization of commutative rings whose artinian modules are RD-injective (Theorem~\ref{T:gloa}). We begin by proving these two theorems for local rings . We show that
the family of proper principal ideals of a local ring $R$ is directed
if and only if $R/P$ is RD-injective,
where $P$ is the maximal ideal of $R$. In this case every noetherian
$R$-module is RD-coflat and each artinian $R$-module is RD-flat. When
$R$ satisfies the ascending chain condition
on principal ideals, then $R$ is a valuation ring if and only if $R/P$ is
RD-injective. However this is not generally true and we give examples of
local rings $R$ which are not valuation rings and such that $R/P$ is
RD-injective. 

All rings in this paper are associative with unity, and all modules are
unital. A left $R$-module is said to be \textbf{cyclically
  presented} if it is of the form $R/Rr$, where $r\in R$. We say that
a left module is \textbf{uniserial} if its set of
submodules is totally ordered by inclusion. Recall that a commutative
ring $R$ is a \textbf{valuation ring} if it is uniserial as
$R$-module and that $R$ is \textbf{arithmetic} if $R_P$ is a
valuation ring for every maximal ideal $P$.

An exact sequence of left $R$-modules $0\rightarrow F\rightarrow
E\rightarrow G\rightarrow 0$ is \textbf{pure-exact} if it remains exact when
tensoring it with any right $R$-module. In this case we say that $F$ is a
\textbf{pure submodule} of $E.$
When $rE\cap F = rF$ for every $r\in R,$ we say that $F$ is an
\textbf{RD-submodule} of $E$(relatively divisible) and that the sequence
is \textbf{RD-exact.}

An $R$-module $F$ is \textbf{pure-injective} (respectively
\textbf{RD-injective}) if for every pure (respectively RD-) exact
sequence
$0\rightarrow N\rightarrow M\rightarrow L\rightarrow 0$ of
$R$-modules, the following sequence
 $0\rightarrow\mathrm{Hom}_R(L,F)
\rightarrow\mathrm{Hom}_R(M,F)\rightarrow\mathrm{Hom}_R(N,F)\rightarrow
0$ is exact.

An $R$-module $F$ is \textbf{pure-projective} (respectively
\textbf{RD-projective}) if for every pure (respectively RD-) exact
sequence
$0\rightarrow N\rightarrow M\rightarrow L\rightarrow 0$ of
$R$-modules, the following sequence
 $0\rightarrow\mathrm{Hom}_R(F,N)
\rightarrow\mathrm{Hom}_R(F,M)\rightarrow\mathrm{Hom}_R(F,L)\rightarrow
0$ is exact.

\section{RD-flatness and RD-coflatness.}
\label{S:flat}
We begin with some preliminary results. As in \cite{FuSa01} we set
$M^{\flat}=\mathrm{Hom}_{\mathbb{Z}}(M,\mathbb{Q}/\mathbb{Z})$ the
\textbf{character} module of $M$.

As in \cite{Ste67} we say that a right module $M$ is {\bf RD-flat} if
for every RD-exact sequence of left modules 
\(0\rightarrow H\rightarrow F\rightarrow
L\rightarrow 0\), the sequence
 \[0\rightarrow M\otimes_R H\rightarrow
M\otimes_R F\rightarrow M\otimes_R L\rightarrow 0\ \mathrm{is\ exact}.\]
 The next
proposition and its proof is similary to that we know for flat
modules. In particular, to prove that $(1)\Rightarrow (5)$, we do a similar
proof as in \cite[Theorem I.1.2]{Laz69}. 
\begin{proposition}
\label{P:flat}
Let $R$ be a ring and $M$ a right module. Then the following
assertions are equivalent.
\begin{enumerate}
\item $M$ is RD-flat
\item $M^{\flat}$ is RD-injective
\item $M\cong P/Q$, where $P$ is RD-projective and $Q$ is a pure
  submodule of $P$
\item Every RD-exact sequence \(0\rightarrow Q\rightarrow P\rightarrow
M\rightarrow 0\) is pure exact
\item $M$ is direct limit of finite direct sums of cyclically
  presented modules.
\end{enumerate}
When $R$ is commutative and $E$ an injective cogenerator, these
conditions are equi\-valent to the following conditions:
\begin{itemize}
\item[(2')] $\mathrm{Hom}_R(M,E)$ is RD-injective.
\item[(6)] $M_P$ is RD-flat over $R_P$ for each maximal ideal $P$.
\end{itemize}
Moreover, each direct limit of right RD-flat modules is RD-flat.
\end{proposition}

 Some examples of commutative artinian rings $R$,
such that $R$ is not an RD-injective module are known:
\cite[Example~1]{BrCo75} and \cite[Example~1]{Pun92}. From Proposition~\ref{P:flat} we deduce the following:

\begin{proposition} \label{P:artring}
Let $R$ be an artinian commutative ring. Then $R$ is an RD-injective module if and only if
$R$ is a quasi-Frobenius ring.
\end{proposition}
\begin{proof}
$R$ is a finite product of local rings. So, we may assume that $R$ is local. Let $E$ be the injective hull of its residue field.
We know that $R\cong\mathrm{Hom}_R(E,E)$. By Proposition~\ref{P:flat} $R$ is an RD-injective module if and only if $E$ is RD-flat.
But it is also known that  $E$ is finitely presented. It follows that $E$ is RD-projective because a module is RD-projective
if and only it is RD-flat and pure-projective. By \cite[Corollary 2]{War69} $E$ is a direct sum of cyclically presented modules. Hence $E\cong R$
 since $E$ is indecomposable and faithful. We get that $R$ is quasi-Frobenius.
\end{proof}

\bigskip
We say that a right $R$-module $M$ is \textbf{RD-coflat} if every RD-exact
sequence
\(0\rightarrow M\rightarrow P\rightarrow
Q\rightarrow 0\) is pure exact. 
\textbf{Pure-essential extension},
  \textbf{RD-essential extension}, \textbf{pure-injective hull} and
  \textbf{RD-injective hull} are defined as in \cite{War69} or
  \cite[Chapter XIII]{FuSa01}. Observe that the RD-injective hull and the pure-injective hull coincide for any
   RD-coflat module. So one gets a partial answer to \cite[Problem~47]{FuSa01}. 
\begin{proposition}
\label{P:cflrd} Let $R$ be a ring. Then:
\begin{enumerate}
\item Pure submodules of RD-coflat modules are RD-coflat. Direct products and direct sums of RD-coflat right
$R$-modules are RD-coflat.
\item for any right module $M$ the following conditions are equivalent:
\begin{itemize}
\item $M$ is RD-coflat.
\item For each RD-exact sequence $0\rightarrow K\overset{u}\rightarrow
  L\overset{p}\rightarrow F\rightarrow 0$ of right modules, where $F$ is pure-projective, the following sequence is exact:
\[0\rightarrow\mathrm{Hom}_R(F,M)\rightarrow\mathrm{Hom}_R(L,M)\rightarrow\mathrm{Hom}_R(K,M)\rightarrow
0.\]
\end{itemize}
\item A right $R$-module $M$ is RD-injective if and only if it is
  RD-coflat and satisfies the following condition: for any family
  of subgroups $(N_i)_{i\in I}$, where $N_i=(Mr_i:s_i)$, for some
  elements $r_i$ and $s_i$ of $R$, and any family $(x_i)_{i\in I}$ of
  elements of $M$, if the sets $x_i+N_i$ have the finite intersection
  property, then their total intersection is non-empty.

\end{enumerate}
\end{proposition}
\begin{proof}
$(1).$ Let $E$ be a right RD-coflat module, $F$ a right module and $M$ a pure
submodule of $E$ which is an RD-submodule of $F$. Let $H$ be the module
defined by the following pushout diagram:
\[\begin{matrix}M&\rightarrow&E\\
\downarrow&{}&\downarrow\\
F&\rightarrow&H\end{matrix}\]
It is easy to prove that $E$ is an RD-submodule of $H$. It successively
follows that $E$ is a pure submodule of $H$, $M$ a pure submodule of
$H$ and $M$ a pure submodule of $F$.

Let $(M_i)_{i\in I}$ be a family of RD-coflat modules
and for each $i\in I$ let $E_i$ be the RD-injective hull of
$M_i$. Then $E=\Pi_{i\in I}E_i$ is RD-injective and $M=\Pi_{i\in
  I}M_i$ is a pure submodule of $E$. It follows that $M$ is
RD-coflat. Whence $\oplus_{i\in I}M_i$ is also RD-coflat since it is a
pure submodule of $M$.

$(2).$ Assume that $M$ is RD-coflat. Let $M\overset{v}\rightarrow\widehat{M}$ be the
RD-injective hull of $M$ and $f:K\rightarrow M$ be a
homormorphism. There exists a morphism $g:L\rightarrow\widehat{M}$
such that $v\circ f=g\circ u$. Then, if
$q:\widehat{M}\rightarrow\widehat{M}/M$ is the natural map, there
exists a morphism $h:F\rightarrow\widehat{M}/M$ such that $h\circ
p=q\circ g$. Since $v$ is a pure monomorphism, there exists a morphism
$j:F\rightarrow\widehat{M}$ such that $h=q\circ j$. It is trivial to
verify that $q\circ (g-j\circ p)=0$. It follows that there exists a
morphism $l:L\rightarrow M$ such that $v\circ l=g-j\circ p$. We easily
get that $v\circ f=v\circ l\circ u$. Since $v$ is a monomorphism, thus
$f=l\circ u$.

Conversely,  let $x_1,\dots,x_n\in\widehat{M}$, and
$\{a_{i,j}\mid 1\leq i\leq n,\ 1\leq j\leq m\}$ be a family of
elements of $R$ such that $\Sigma_{i=1}^{i=n}x_ia_{i,j}=y_j\in M$, for
each $j,\ 1\leq j\leq m$. Let $F$ be a finitely
presented module generated by $e_1,\dots,e_n$ with the relations
$\Sigma_{i=1}^{i=n}e_ia_{i,j}=0$ for each $j,\ 1\leq j\leq m$ and
$h:F\rightarrow\widehat{M}/M$ be the morphism defined by
$h(e_i)=q(x_i),\ \forall i,\ 1\leq i\leq n$. We
consider the following pullback diagram:
\[\begin{matrix}L&\overset{p}\rightarrow&F\\
\downarrow&{}&\downarrow\\
\widehat{M}&\overset{q}\rightarrow&\widehat{M}/M \end{matrix}\]
Let $g$ be the left vertical map. Then $M$ is isomorphic to $\ker
p$. It follows that $M$ is an 
RD-submodule of $L$. Therefore there exists a morphism
$l:L\rightarrow M$ such that $l(x)=x,$ for each $x\in M$. Let
$z_1,\dots,z_n$ be elements of $L$ such that $g(z_i)=x_i$ and
$p(z_i)=e_i$ for each $i,\ 1\leq i\leq n$. Then we get that
$\Sigma_{i=1}^{i=n}l(z_i)a_{i,j}=y_j$ for each $j,\ 1\leq j\leq
m$. Hence $M$ is a pure submodule of $\widehat{M}$.

$(3)$. We do the same proof as in \cite[Theorem~4]{War69}.
\end{proof}
 
\bigskip
From these previous propositions and \cite[Theorem~2.5]{PPR99} we
deduce the follo\-wing:
\begin{theorem}
\label{T:flat}
Let $R$ be a ring. The following assertions are equivalent:
\begin{enumerate}
\item Every right $R$-module is RD-flat
\item Every left $R$-module is RD-coflat.
\item Every left pure-injective module is RD-injective
\item Every right pure-projective module is RD-projective
\item Every right finitely presented module is a summand of a direct
  sum of cyclically presented modules.
\item Every RD-exact sequence of right modules is pure-exact.
\item The left-right symmetry of $(1)-(6)$.
\end{enumerate}
\end{theorem}

As in \cite{PPR99} we will say that $R$ is an {\bf RD-ring} if it satisfies
the equivalent conditions of Theorem~\ref{T:flat}. By
\cite[Proposition~4.5]{PPR99} a commutative ring $R$ is RD if and
  only if it is an arithmetic ring.

\section{Rings whose RD-injective right modules are $\Sigma$-RD-injective.}
\label{S:hull}

In this section we will prove the following Theorem~\ref{T:cofinj}. In
the sequel, for every right $R$-module $M$, we denote its RD-injective hull by
$\widehat{M}$.

As in
\cite[Observation~3(2)]{HZ00}, if $N$ is a subgroup of a right module
$M$ over a ring $R$, we say
that $N$ is \textbf{a finite matrix subgroup} if it is the kernel of the
following map $M\rightarrow M\otimes_R X$, defined by $m\rightarrow
m\otimes x$, where $X$ is a finitely presented left module, $x\in X$
and $m\in M$. For instance, for any $r,s\in R$, $(Mr:s)$ is a finite matrix subgroup: it is the kernel of the following map: $M\rightarrow
M\otimes_R R/Rr$ defined by $m\rightarrow m\otimes s$. If $s$ is a
unit then $(Mr:s)=Mr$.

Recall that a right $R$-module $M$ is $\Sigma$-\textbf{pure-injective}
if $M^{(I)}$ is pure-injective for each index set $I$. A ring $R$ is
said to be \textbf{right perfect} if every flat right
$R$-module is projective.

Recall that a right module $M$ is \textbf{fp-injective} (or
\textbf{absolutely pure}) if it is a pure submodule of each
overmodule.
\begin{theorem} 
\label{T:cofinj}
Let $R$ be a ring. Consider the following conditions:
\begin{enumerate}
\item Every RD-coflat right $R$-module is RD-injective.
\item Each direct sum of right RD-injective modules is RD-injective.
\item Every RD-injective right $R$-module is $\Sigma$-pure-injective.
\item Every RD-injective right $R$-module is a direct sum of
  indecomposable submo\-dules.
\item $R$ is right artinian and RD-injective hulls of finitely
  generated right mo\-dules are finite direct sums of indecomposable
  submo\-dules.  
\item $R$ is a
  finite product of artinian valuation rings or finite rings.
\end{enumerate}
Then:
\begin{itemize} 
\item Conditions $\mathrm{(1),\ (2),\ (3)}$ and $\mathrm{(4)}$ are equivalent and imply condition $\mathrm{(5)}$.
\item When $R$ is commutative, the six conditions are equivalent. 
\end{itemize}
\end{theorem}
\begin{proof}
$(1)\Rightarrow (2)$ follows from Proposition~\ref{P:cflrd}. 

$(2)\Rightarrow (3)$ is obvious.

$(3)\Rightarrow (1)$. Let $M$ be a right RD-coflat module. Then $\widehat{M}$ is
$\Sigma$-pure-injective. We conclude by \cite[Corollary~8]{HZ00}.

$(4)\Leftrightarrow (3)$ follows from \cite[Proposition~9 and Theorem~10]{HZ00}. 

$(2)\Rightarrow (5)$. Each direct sum of right injective modules is fp-injective and RD-coflat, hence injective. We deduce that $R$ is right
noetherian. Since $\widehat{R}$ is $\Sigma$-pure-injective then it
satisfies the descending chain condition on its finite matrix
subgroups by \cite[Observation~3(4) and Theorem~6]{HZ00}. Consequently the family
$(\widehat{R}r)_{r\in R}$ satisfies the descending chain condition. Since
$Rr=\widehat{R}r\cap R$ for each $r\in R$, it follows that $R$ verifies this
descending chain condition on principal left ideals. Hence $R$ is also
right perfect by \cite[Th\'eor\`eme~6.2.5]{Ren75}. So $R$ is right artinian by
\cite[Proposition~6.2.10]{Ren75}.

Let $M$ be a finitely generated right $R$-module. Since $(2)\Rightarrow
(4)$, $\widehat{M}=\oplus_{i\in I}M_i$ where $M_i$ is indecomposable for
each $i\in I$. There is a finite subset $J$ of $I$ such that
 $M$ is relatively divisible in
$\oplus_{i\in J}M_i$. Since $M\cap\oplus_{i\in I\setminus J}M_i=0$ and
$\oplus_{i\in J}M_i\cong \widehat{M}/\oplus_{i\in I\setminus J}M_i$,
we conclude that $\oplus_{i\in I\setminus J}M_i=0$ and $I=J$. 
\end{proof}
 
\bigskip
 To prove the second assertion of this theorem some
preliminary results are needed.
\begin{proposition} 
\label{P:finite}
Let $R$ be a finite ring. Then every RD-coflat
  right (or left) module is RD-injective.
\end{proposition} 
\begin{proof}
For every right $R$-module $M$ the family of
subgroups which are finite intersections of subgroups of the form
$(Mr:s)$, where $r,s\in R$, is finite. We conclude by
Proposition~\ref{P:cflrd}.
\end{proof}

\bigskip
Now we consider a commutative local ring $R$ of maximal ideal $P$, satisfying the ascending
chain condition on principal ideals . If $A$ is a proper ideal we
denote $A^*$ the set of principal ideals contained in $A$, $\mathcal{A}^*$ the subset of maximal elements of $A^*$. We put 
$S=R/P$ and $E=\mathrm{E}_R(S)$. For any
$r^*\in A^*$ let $\phi_{r^*}:R/r^*\rightarrow R/A$ be the homomorphism
defined by $\phi_{r^*}(1+r^*)=1+A$. Let $\psi:\oplus_{r^*\in \mathcal{A}^*}R/r^*\rightarrow R/A$
be the homomorphism defined by the family $(\phi_{r^*})_{r^*\in
  \mathcal{A}^*}$. We put $\Psi=\mathrm{Hom}_R(\psi,E)$ and
\(E[A]=\{e\in E\mid A\subseteq (0:e)\}.\) 
Recall that $\mathrm{Hom}_R(R/A,E)\cong E[A]$. It is easy to
check that $\Psi$ is the diagonal homomorphism induced by the
inclusion maps $E[A]\rightarrow E[r^*]$. The following lemma
holds.
\begin{lemma} Then:
\label{L:simp}
\begin{enumerate}
\item $\ker\psi$ is an RD-submodule of $\oplus_{r^*\in \mathcal{A}^*}
R/r^*$. 
\item $\Psi(E[A])$ is an RD-submodule of $\prod_{r^*\in \mathcal{A}^*}E[r^*]$.
\end{enumerate}
\end{lemma}
\begin{proof} $(1)$ Let  
$f:R/s^*\rightarrow R/A$ be a non-zero homomorphism where $s^*$ is a proper principal ideal of $R$. We have
$f(1+s^*)=a+A$ where $a\in R\setminus A$. Then $as^*\subseteq A$. Let $t^*\in\mathcal{A}^*$ such that $as^*\subseteq t^*$. Let
$g:R/s^*\rightarrow \oplus_{r^*\in\mathcal{A}^*}R/r^*$ be the composition of
$\alpha:R/s^*\rightarrow R/t^*$, defined by $\alpha(1+s^*)=a+t^*$,
with the inclusion map $R/t^*\rightarrow\oplus_{r^*\in\mathcal{A}^*}R/r^*$. Then $f=\phi\circ g$. We conclude by
\cite[Proposition 2]{War69}. 

$(2)$ is a consequence of the fact that $\mathrm{Hom}_R(M/N,E)$ is isomorphic to an RD-submodule of $\mathrm{Hom}_R(M,E)$
if $N$ is an RD-submodule of $M$.
\end{proof}

\bigskip
Since its maximal ideal is T-nilpotent, it is easy to prove that every
commutative local perfect ring satisfies the ascending chain
condition on principal ideals.

\begin{lemma}
\label{L:hull}
Let $R$ be a commutative local perfect ring and $A$ a proper ideal. 
Then
$\widehat{E[A]}\cong\prod_{r^*\in\mathcal{A}^*}E[r^*]$.  
\end{lemma}
\begin{proof} We put $G(A)=\prod_{r^*\in\mathcal{A}^*}E[r^*]$. By \cite[Lemma XIII.1.2]{FuSa01} $E[r^*]$ is RD-injective, so $G(A)$ is RD-injective too. If $N$
is a submodule of $G(A)$ such that
$N\cap\Psi(E[A])=0$ and the image of $\Psi(E[A])$ relatively
divisible in $G(A)/N$, then 
$\forall 0\not=y\in N$, $Ry\cap\Psi(E[A])=0$ and the image of $\Psi(E[A])$
is relatively divisible in
$G(A)/Ry$. Consequently it is sufficient
to prove that for every $y\in G(A)$ such that
$Ry\cap\Psi(E[A])=0$, the image of $\Psi(E[A])$ is not relatively
divisible in $G(A)/Ry$. Let
$0\not=y=(y_{r^*})_{r^*\in\mathcal{A}^*}$ such that
$Ry\cap\Psi(E[A])=0$. Let $F$ be the submodule of $E$ generated by
$\{y_{r^*}\mid r^*\in \mathcal{A}^*\}$, $B=\mathrm{ann}(F)$ and $t\in R$ such
that $t+B$ generates a minimal non-zero submodule of $R/B$. It follows
that $tF=S$. Hence there exists $s^*\in\mathcal{A}^*$ such that
$ty_{s^*}\not=0$. We set $e=ty_{s^*}$ and $z_{s^*}=0$. If
$s^*\not=r^*\in \mathcal{A}^*$ and $s^*=Rs$ then $s\notin r^*$. It
follows that $(r^*:s)\subseteq P\subseteq (0:e-ty_{r^*})$. Since
$E[r^*]$ is injective over $R/r^*$ there exists $z_{r^*}\in E[r^*]$
such that $sz_{r^*}=e-ty_{r^*}$. We put $z=(z_{r^*})_{r^*\in
  \mathcal{A}^*}$. Then we get the following equality:
$sz+ty=\Psi(e)$. Since $Ry\cap\Psi(E[A])=0$ we have $sz\not=0$. If
$e=sx$ for some $x\in E[A]$ then $e=0$ since $s\in A$. It follows that
the image of $\Psi(E[A])$ is not relatively divisible in
$G(A)/Ry$.
\end{proof}
 
\begin{remark} \textnormal{We know that $E[A]$ is indecomposable
    and it is a module of finite length if $R$ is a commutative local artinian
    ring. If $A$ is not principal then $\widehat{E[A]}$ is
    decomposable, and if $\mathcal{A}^*$ is not finite then
    $\widehat{E[A]}$ is neither artinian nor noetherian.}
\end{remark}
\begin{remark} \textnormal{Let $\alpha:F\rightarrow M$ be an epimorphism
    of right modules where $F$ is RD-projective. Let
    $K=\ker\alpha$. We say that 
    $\alpha:F\rightarrow M$ is an {\bf RD-projective cover} of $M$ if
    $K$ is an RD-pure submodule of $F$ and $F$ is the only submodule
    $N$ of $F$ which verifies $K+N=F$ and $K\cap N$ is relatively
    divisible in $N$. In a similar way we define the {\bf pure-projective cover} of $M$. If $R$ is a commutative local perfect ring,
    then Lemma~\ref{L:hull} implies that, for each ideal $A$,
    $\psi:\oplus_{r^*\in\mathcal{A}^*}R/r^*\rightarrow R/A$ is an
    RD-projective cover of $R/A$.} 
\end{remark}

Now we prove the last assertion of Theorem~\ref{T:cofinj}.

\begin{proof}[Proof of Theorem~\ref{T:cofinj}.]
$(5)\Rightarrow (6)$. $R$ is a finite product of local rings. We may
assume that $R$ is local. There exists a finite family of
indecomposable RD-injective modules $(S_i)_{1\leq i\leq n}$ such that
$\widehat{S}=\oplus_{i=1}^{i=n}S_i$. By \cite[Proposition~9(3)]{HZ00} the
endomorphism ring of each indecompo\-sable pure-injective module is 
local. For every $r^*\in\mathcal{P}^*$, $E[r^*]$ is a summand of
$\widehat{S}$ by Lemma~\ref{L:hull}. Then there exists $i,$ $1\leq
i\leq n$ such that $E[r^*]\cong S_i$ by Krull-Schmidt Theorem. Whence
$\mathcal{P}^*$ is a finite set. This implies that the set
$\mathcal{P'}^*$ of maximal principal proper ideals of $R/P^2$ is
finite. Clearly $\mathcal{P'}^*$ is the set of vector lines of the
vector space $P/P^2$ over $R/P$. This last set is finite if and only
if $P$ is principal or $R/P$ is a finite field. 

$(6)\Rightarrow (1)$. We may assume that $R$ is local. If $R$ is an
artinian valuation ring then each
$R$-module is RD-injective by \cite[Theorem~4.3]{Gri70}. We conclude by
Proposition~\ref{P:finite} if $R$ is finite.
\end{proof}

\section{Rings whose RD-flat modules are RD-projective.}
\label{S:perf}

By using a result of Crawley-Boevey in \cite{CrBo94} we get the next theorem. If $R$ is a ring, we denote by $J$ its Jacobson radical and by $\mathcal{RD}$ the family of proper right ideals $A$ such that $R/A$ is isomorphic to an indecomposable summand of a cyclically presented right module.

\begin{theorem} \label{T:perf}
Let $R$ be a ring. The following conditions are equivalent:
\begin{enumerate}
\item each RD-flat right $R$-module is RD-projective
\item $R$ is right perfect and $\mathcal{RD}$ satisfies the ascending chain condition.
\end{enumerate}
\end{theorem}
\begin{proof}
The condition $(1)$ implies that each flat right module is projective. Hence $R$ is right perfect.
Let $\mathcal{A}$ be the category of RD-flat right $R$-modules. By Proposition~\ref{P:flat}.5 each object of $\mathcal{A}$ is direct limit of finitely presented objects of $\mathcal{A}$. So, $\mathcal{A}$ is an accessible additive category (or a locally finitely presented additive category). Then each object of $\mathcal{A}$ is pure-projective. By \cite[(3.2) Theorem]{CrBo94} we get:
\begin{itemize}
\item[(a)] Every finitely presented object of $\mathcal{A}$ is a finite direct sum of indecomposable objects, each with local endomorphism ring.
\item[(b)] Given a sequence
\begin{equation} \label{E:seq}
X_1\stackrel{f_1}{\rightarrow}X_2\stackrel{f_2}{\rightarrow}X_3\stackrel{f_3}{\rightarrow}X_4\stackrel{f_4}{\rightarrow}\dots \tag*{(S)}
\end{equation}
of non-isomorphism between finitely presented indecomposable objects of $\mathcal{A}$, the composition $f_N\circ\dots\circ f_2\circ f_1$ is zero for $N$ sufficiently large.
\end{itemize}
Suppose that $\mathcal{RD}$ contains a strictly ascending chain $(A_n)_{n\in\mathbb{N}^*}$. Let $X_n=R/A_n$ and let $f_n:X_n\rightarrow X_{n+1}$ be the natural map, $\forall n\in\mathbb{N}^*$. Then we get a contradiction since $f_n\circ\dots\circ f_2\circ f_1\ne 0,\ \forall n\in\mathbb{N}^*$.

Conversely\footnote{The proof of $(2)\Rightarrow (1)$ is valid if we assume that $R$ is local. We shall give a complete proof in the last section} we must prove that $\mathcal{A}$ satifies the above conditions (a) and (b).
If we show that $\mathrm{End}_R(X)$ is local for each indecomposable summand of a cyclically presented right module then the condition (a) is an immediate consequence of Krull-Schmidt Theorem. Let $x$ be a generator of $X$ and $f\in\mathrm{End}_R(X)$. First we show that $f$ is an isomorphism if $f$ is onto. For each integer $n\geq 1$ we put $Y_n=\mathrm{Ker}\ f^n$ and $A_n=\{r\in R\mid xr\in Y_n\}$. Since $X/Y_n\cong X$, $A_n\in\mathcal{RD}$ for each integer $n\geq 1$. Thus there exists an integer $N$ such that $A_n=A_{n+1},\ \forall n\geq N$. Consequently  $Y_n=Y_{n+1},\ \forall n\geq N$. We easily deduce that $Y_n\cap f^n(X)=0$ if $n\geq N$. Since f is surjective we get that $f$ is an isomorphism. Now suppose that $f$ is not surjective. Then $x-f(x)$ is a generator of $X$ because\footnote{this is true if $R$ is local} $f(x)\in XJ$. It follows that $\mathbf{1}_X-f$ is onto and consequently it is an isomorphism. So, $\mathrm{End}_R(X)$ is local.

Now we consider the sequence~\ref{E:seq}. Let $x_n$ be a generator of $X_n$ and let $t_n\in R$ such that $f_n(x_n)=x_{n+1}t_n$, for each integer $n\geq 1$. Then the set of integers $n$ for which $t_n\in J$ is infinite\footnote{this is true if $R$ is local}. Otherwise there exists an integer $m$ such that $t_n\notin J,\ \forall n\geq m$. By induction on $n$ we define $y_n\in X_n$ for each integer $n\geq m$ in the following way: $y_m=x_m$ and $y_{n+1}=f_n(y_n)$. Let $A_n$ be the annihilator of $y_n$. Then $y_n$ is a generator of $X_n$ for each integer $n\geq m$ and we have $A_n\subset A_{n+1}$ since $f_n$ is not injective. This contradicts that $\mathcal{RD}$ satisfies the ascending chain condition.  Now, from the T-nilpotence of $J$ we deduce that there exists an integer $N$ such that $t_N\dots t_2t_1=0$, whence  $f_N\circ\dots\circ f_2\circ f_1$ is zero. The proof is now complete.
\end{proof}

\medskip
We easily deduce the two following corollaries.
\begin{corollary} \label{C:local}
Let $R$ be a local ring. Then $R$ is right perfect if and only if each RD-flat right module is RD-projective.
\end{corollary}
\begin{proof}
 In this case, $\mathcal{RD}$ is the family of right principal ideals. From the right T-nilpotence of $J$ we deduce that this family satisfies the ascending chain condition. \end{proof}

\begin{corollary} \label{C:perf}
Let $R$ be a commutative ring. Then $R$ is perfect if and only if each RD-flat module is RD-projective.
\end{corollary}
\begin{proof}
If $R$ is perfect then $R$ is a finite product of local rings.
\end{proof}

\begin{remark} \textnormal{Observe that each right artinian ring satisfies the equivalent conditions of Theorem~\ref{T:perf}. From Theorem~\ref{T:cofinj} it follows that the condition "each RD-coflat right $R$-module is RD-injective" implies the condition "each RD-flat right $R$-module is RD-projective". But the converse doesn't hold.}

\textnormal{If $R$ is a local right perfect ring which is not left perfect then each right RD-flat module is RD-projective but it is not true for left modules. An example of a such ring is given in \cite[Exemple 3, p.132]{Ren75}.}
\end{remark}

\section{RD-injective artinian modules.}
\label{S:inje}
The aim of this section is to characterize commutative rings for
which every artinian module is RD-injective. Throughout this section
 all rings are commutative. Recall that a family $\mathcal{F}$ of
 submodules of an $R$-module $M$ is {\bf directed} if, $\forall (U,V)\in
 \mathcal{F}^2$ there exists $W\in\mathcal{F}$ such that $U+V\subseteq
 W$. We will prove the two following theorems.

\begin{theorem} 
\label{T:glo}
Let $R$ be a ring. Then the following conditions are
  equivalent:
\begin{enumerate}
\item Every simple $R$-module is RD-flat.
\item Every simple $R$-module is RD-injective.
\item Every $R$-module of finite length is RD-flat.
\item Every $R$-module of finite length is RD-injective.
\item Every artinian $R$-module is RD-flat.
\item Every noetherian $R$-module is RD-coflat.
\item For each maximal ideal $P$ the
  family of proper principal ideals of $R_P$ is directed.
\end{enumerate}
Moreover, if $R$ satisfies these conditions, then every artinian module is a
(finite) direct sum of uniserial modules and every noetherian module
is a direct sum of $2$-generated submodules. However the converse is not true.
\end{theorem}
\begin{proof}
It is obvious that $(5)\Rightarrow (3)\Rightarrow (1)$ and $(6)\Rightarrow (4)\Rightarrow (2)$.

$(1)\Leftrightarrow (2)$. We set $E=\prod_{P\in\mathrm{Max(R)}}\mathrm{E}_R(R/P)$ where
$\mathrm{Max}(R)$ is the set of maximal ideals of $R$. Then $E$ is an injective
cogenerator and for each simple module $S$, $S\cong \mathrm{Hom}_R(S,E)$. We
conclude by Proposition~\ref{P:flat}.

$(3)\Rightarrow (5)$ since each artinian
module is direct limit of modules of finite length.
\end{proof}
 
\bigskip
 Let $R$ be a ring. For every maximal ideal $P$ we denote by
$L(P)$ the sum of all submodules of finite length of $\mathrm{E}(R/P)$ and we
set $J(P)=\cap_{n\in\mathbb{N}}(PR_P)^n$.  

\begin{theorem}
\label{T:gloa}
Let $R$ be a ring. Then the following conditions are
  equivalent:
\begin{enumerate}
\item Every artinian $R$-module is RD-injective.
\item Every noetherian $R$-module is RD-flat.
\item $R/P$ and $L(P)$ are RD-injective for every maximal ideal $P$.
\item For each maximal ideal $P$,  the
  family of proper principal ideals of $R_P$ and the family of
  principal ideals of $R_P$ contained in $J(P)$ are directed.
\end{enumerate}
Moreover, if $P=P^2$  for each maximal ideal $P$, then these
conditions are equivalent to conditions of
Theorem~\ref{T:glo}. In this case, every noetherian module and every
artinian module is semi-simple.
\end{theorem}

To prove this theorem and complete the proof of Theorem~\ref{T:glo} we
study our problem in the local case.

In the sequel we assume that $R$ be a local ring. We denote $P$ its
maximal ideal, $E$ the injective hull of $R/P$ and
$J=\cap_{n\in\mathbb{N}}P^n$. Let $A$ be a proper ideal of $R$ and $A^*$
the set of non-zero principal ideals contained in $A$. Let $\phi:\oplus_{r^*\in A^*}
R/r^*\rightarrow R/A$ be the homomorphism defined by the family
$(\phi_{r^*})_{r^*\in A^*}$. The two
following lemmas hold. The first is similar to Lemma~\ref{L:simp}.
\begin{lemma} 
\label{L:sim}
Then $\ker\phi$ is an RD-submodule of $\oplus_{r^*\in A^*}R/r^*$.
\end{lemma}

\begin{lemma}
\label{L:main}
Let $R$ be a local ring, and $A$ a proper ideal. Then:
\begin{enumerate}
\item $E[A]$ is RD-injective if and only if $R/A$ is RD-flat
\item $R/A$ is RD-flat if and only if $A^*$ is directed.
\item If $A$ is finitely generated then $E[A]$ is RD-injective if and
  only if $A$ is principal.
\end{enumerate}
\end{lemma}
\begin{proof}
(1) is an immediate consequence of Proposition~\ref{P:flat}.

Assume that $R/A$ is RD-flat. By Proposition~\ref{P:flat}
$\ker\phi$ is a pure submodule of $\oplus_{r^*\in
  A^*}R/r^*$. Let $I\subseteq A$ be a finitely generated ideal and
$f:R/I\rightarrow R/A$ be the homomorphism defined by $f(1+I)=1+A$. By
\cite[Proposition 3]{War69}, for some morphism $g:R/I\rightarrow
\oplus_{r^*\in A^*}R/r^*$, $f=\phi\circ g$. Then $g(1+I)=\sum_{r^*\in
  F}a_{r^*}+r^*$ where $F$ is a
finite family of $A^*$. We deduce that $1+A=\sum_{r^*\in F}a_{r^*}(1+A)$. Since
$R$ is local there exists $r^*\in F$ such that
$a_{r^*}$ is a unit of $R$. It follows that $I=(0:1+I)\subseteq
(0:a_{r^*}+r^*)=r^*$.

Conversely if $R/A=\varinjlim_{r^*\in A^*}R/r^*$ then $R/A$ is RD-flat
by Proposition~\ref{P:flat}.
\end{proof}

\bigskip
The next theorem characterizes local rings for which every simple
module is RD-injective (and RD-flat). It is an immediate consequence of
Lemma~\ref{L:main}.  
\begin{theorem} 
\label{T:main}
Let $R$ be a local ring. 
Then $R/P$ is RD-flat if and only if $P^*$ is directed. 
\end{theorem}

 The following corollaries will be useful to provide examples
of rings for which every simple module is RD-flat.
\begin{corollary}
\label{C:cons}
Let $R$ be a local ring. Consider the two following conditions:
\begin{enumerate}
\item $R/P$ is RD-flat.
\item There exists an ideal $I\subset P$ such
  that $R/I$ is a valuation ring and for each $r\in R\setminus I$,
  $I\subset Rr$.
\end{enumerate}
Then $(2)\Rightarrow (1)$. When $P$ is finitely generated the two
conditions are equivalent.
\end{corollary}
\begin{proof}
First we assume that $R/P$ is RD-flat and $P$ is
finitely generated. By
Theorem~\ref{T:main} $P$ is principal. We put $P=Rp$ and 
$I=J$. Let $a\in P\setminus J$. There
exists $n\in\mathbb{N}$ such that $a\notin Rp^n$. We may assume that $n$ is
minimal. Then $a=bp^{n-1}$ and since $a\notin Rp^n$, $b$ is a unit of
$R$. Hence $I\subset Rp^{n-1}=Ra$ and $R/I$ is a
valuation ring.
Conversely, if $a$ and $b$ are elements of $P\setminus I$ it is easy
to prove that $Ra\subseteq Rb$ or $Rb\subseteq Ra$. By
Theorem~\ref{T:main} $R/P$ is RD-flat.
\end{proof}

\begin{corollary}
Let $R$ be a local ring. Assume that $R$ satisfies the ascending chain
condition on principal ideals.
Then $R$ is a 
valuation ring if and only if $R/P$ is RD-injective.
\end{corollary}
\begin{proof}
Assume that $R/P$ is RD-injective. We claim that
there exists only one maximal principal ideal $pR$ and that $pR=P$. If
there exist at least two maximal principal ideals $p^*$ and $q^*$ then
$p^*+q^*$ is strictly contained in a proper principal ideal. Consequently we
obtain a contradiction. If $J\ni a\not=0$
then, by induction, we get a strictly ascending chain of ideals
$(Rc_n)_{n\in\mathbb{N}^*}$ such that $a=pc_1$ and
$c_n=pc_{n+1},\forall n\in\mathbb{N}^*$. Hence $J=0$. We complete
the proof by using Corollary~\ref{C:cons}.
\end{proof}

\bigskip
We give two examples of local rings $R$ which are not valuation rings
and such that $R/P$ is RD-injective.
\begin{example} 
\label{E:fg} \textnormal{Let $D$ a valuation domain, $Q$ its field
    of fractions, $H =Q^2$.}

 \textnormal{We put $R =
\displaystyle\Bigl\{{\binom{d\ \ h}{0\ \ d}}\mid d\in D, h\in
H\Bigr\}$ the trivial extension of $D$ by $H.$}

 \textnormal{If $I=\displaystyle\Bigl\{{\binom{0\ \ h}{0\ \ 0}}\mid h\in
H\Bigr\}$ then $R/I\cong D$. By using that $H$ is a divisible
$D$-module, it is easy to prove that $I\subset Ra$ for every $a\in
R\setminus I$.}
\end{example}
\begin{example} 
\label{E:nfg} \textnormal{Let $K$ be a field and $T$ the factor ring
 of the polynomial ring $K[X_{p,n}\mid
    p,n\in\mathbb{N}]$ modulo the ideal generated
    by $\{X_{p,n}-X_{p,n+1}X_{p+1,n}\mid
    p,n\in\mathbb{N}\}$.  We
    denote $x_{p,n}$ the image of $X_{p,n}$ in $T$ and $P$ the maximal ideal of
    $T$ generated by $\{x_{p,n}\mid
    p,n\in\mathbb{N}\}$. We put $R=T_P$. Let $p,q,m,n\in\mathbb{N}$ 
    . Clearly $x_{p,n}\in Rx_{p+q,n}\subset Rx_{p+q,n+m}$ and
    $x_{q,m}\in Rx_{p+q,m}\subset Rx_{p+q,n+m}$. If $a\in R$ then
    there exists $b\in T$ and $s\in T\setminus P$ such that
    $a=\displaystyle{\frac{b}{s}}$. Then $b$ is a linear combination
    with coefficients in $K$ of finite
    products of elements of $\{x_{p,n}\mid
    p,n\in\mathbb{N}\}$. From above we deduce that
    $a\in Rx_{p,n}$ for some
    $(p,n)\in\mathbb{N}^2$. It follows from this
    property that every proper finitely generated ideal is contained
    in a proper principal ideal.}

\textnormal{If we prove that $Rx_{p+1,n}$ and $Rx_{p,n+1}$ are not
  comparable for every $(p,n)\in\mathbb{N}^2$ then $R$ doesn't satisfy
  the condition (2) of Corollary~\ref{C:cons}. Else there exists
  $p,n\in\mathbb{N}$ such that $x_{p,n}\notin I$. It follows that
  $x_{p+1,n}$ and $x_{p,n+1}$ are not elements of $I$, whence
  $Rx_{p+1,n}$ and $Rx_{p,n+1}$ are
  comparable. We get a contradiction.}

\textnormal{For every $\ell\in\mathbb{N}$, let
  $T_{\ell}=K[X_{p,\ell-p}\mid 0\leq p\leq\ell]$,  
  $\varphi_{\ell}:T_{\ell}\rightarrow T_{\ell+1}$ be the homomorphism
  defined by
  $\varphi_{\ell}(X_{p,\ell-p})=X_{p+1,\ell-p}X_{p,\ell-p+1}$,
  $\forall p$, $0\leq p\leq\ell$ and $P_{\ell}$ be the ideal of $T_{\ell}$
  generated by $\{X_{p,\ell-p}\mid 0\leq p\leq\ell\}$. Clearly $T$
  is the direct limit of 
  $(T_{\ell},\varphi_{\ell})_{\ell\in\mathbb{N}}$. Now assume that
  $x_{p,n+1}\in Rx_{p+1,n}$. Then there exist $s\in T\setminus
  P$ and $c\in T$ such that $sx_{p,n+1}=cx_{p+1,n}$. There
  exists an integer $\ell\geq p+n+1$ such that $c$ and $s$ are elements of
  $T_{\ell}$. The following formul\ae{} hold:
  $x_{p,n+1}=\prod_{j=0}^{j=\ell-n-p-1}X_{p+j,\ell-p-j}^{\binom{\ell-n-p-1}{j}}$
  and
  $x_{p+1,n}=\prod_{j=0}^{j=\ell-n-p-1}X_{p+j+1,\ell-p-j-1}^{\binom{\ell-n-p-1}{j}}$.\\
  Clearly $x_{p,n+1}$ and $x_{p+1,n}$ are two different monomials of
  $T_{\ell}$. We have $s=k+t$ and $c=a+b$ where $k,a\in K$, $t,b\in
  P_{\ell}$ and $k\not= 0$. We get the following equality:
  $ax_{p+1,n}-kx_{p,n+1}+bx_{p+1,n}-tx_{p,n+1}=0$. It is obvious that the
  degree (in $T_{\ell}$) of $x_{p+1,n}$ and $x_{p,n+1}$ is less than the degree of
  each monomial of $bx_{p+1,n}$ and $tx_{p,n+1}$. It follows that $k=a=0$. Hence we get a
  contradiction. In the same way we prove that $x_{p+1,n}\notin Rx_{p,n+1}$.}
\end{example}

Now we study the RD-injectivity of artinian modules over local rings.
\begin{theorem} 
  \label{T:artin} The following conditions are equivalent for a local ring $R$:
\begin{enumerate}
\item Every artinian $R$-module is RD-injective.
\item $R/P$ and $R/J$ are RD-flat.
\item $R/P$ and $L(P)$ are RD-injective.
\item $P^*$ and $J^*$ are directed.
\end{enumerate}
\end{theorem}
\begin{proof}
Assume that $P$ is not finitely generated and $R/P$ is RD-injective. Let $a\in P$. Then there exists $b\in P\setminus
Ra$. We deduce that $Ra\subset Ra+Rb\subseteq Rc$ for some $c\in
P$. It follows that $a=cd$ for some $d\in P$, whence $P=P^2$. In this
case every artinian module is semi-simple of finite length and $L(P)$
is simple.

 Assume that $P$ is finitely generated and $R/P$ is RD-injective. By
Corollary~\ref{C:cons} and its proof, $P=pR$ for some $p\in P$, and 
$R'=R/J$ is a noetherian valuation ring. Let $M$ be an indecomposable artinian
$R$-module. For every $x\in M$, $Rx$ is a
module of finite length. It follows that there exists
$n\in\mathbb{N}$ such that $p^n\in (0:x)$. We deduce that $M$ is an
$R'$-module. If $U=\mathrm{E}_{R'}(M)$ then $U\cong\mathrm{E}_{R'}(R/P)^m=E[J]^m$ for some
integer $m$. We put $U=U_1\oplus\cdots\oplus U_m$ where $U_j\cong E[J]$ for
each $j$, $1\leq j\leq m$. Thus $M\cap U_j\not=0$ for every $j$,
$1\leq j\leq m$. By \cite[Proposition XII.2.1]{FuSa01} there exists $i$,
$1\leq i\leq m$, such that $M\cap U_i$ is a pure $R'$-submodule of
$M$. It follows that $M\cap U_i$ is a summand of $M$. Hence $M$ is
isomorphic to a submodule of $E[J]$ which is uniserial by \cite[Theorem]{Gil71}. If $M\not=E[J]$ there exists a positive integer $n$
such that $Rp^n=\mathrm{ann}(M)$. In this case we easily deduce that $M$ is
RD-injective. If $M=E[J]$, by Proposition~\ref{P:flat} $M$ is
RD-injective if and only if $R/J$ is RD-flat. Let us observe that
$E[J]=L(P)$. The proof is now complete.
\end{proof}
 
\bigskip
 The following corollary will be useful to provide examples
of rings for which every artinian module is RD-injective.
\begin{corollary}
\label{C:artin} Let $R$ be a local ring. Assume that each artinian
$R$-module is RD-injective and $P$ is finitely generated. Suppose that
$J\not= J^2$. Let $I$ be the inverse image of
$I'=\cap_{n\in\mathbb{N}}(JR_J)^n$ by the natural map $h:R\rightarrow R_J$.
Then $R/I$ is a discrete valuation ring and for every $r\in R\setminus
I$, $I\subset rR$. Moreover, if $J$ is nilpotent, then $R$ is a
discrete valuation ring.
\end{corollary}
\begin{proof}
Let us observe that $J\not=0$ and it is a prime
ideal. By Corollary~\ref{C:cons} $R/J$ is a discrete rank one
valuation domain and $J\subset rR$ for each $r\in R\setminus J$. Let
$q\in J\setminus J^2$ and $r\in J$. Then $Rq+Rr\subseteq Rc$ for some $c\in
J\setminus J^2$. There exist $s\in R\setminus J$ and $t\in R$ such that
$q=sc\hbox{ and }r=tc$. It follows that
$\displaystyle{\frac{r}{1}=\frac{tq}{s}}$. Consequently
$JR_J$ is principal over $R_J$ and $R_J/I'$ is a valuation ring by
Corollary~\ref{C:cons}. We may assume that $I=0$, $I'=0$ and $h$ is a monomorphism. For every $0\not= a\in J$ there exist a unit $u$ and
two integers $m\geq 1$ and $\alpha\geq 0$ such that
$\displaystyle{\frac{a}{1}=\frac{uq^m}{p^{\alpha}}}$. There exists $c\in
J\setminus J^2$ such that $q=cp^{\alpha}$. It follows that
$a=ucq^{m-1}$. If $0\not= b\in J$, in the same way there exist two
integers $n\geq 1$ and $\beta\geq 0$, a unit $v$ and $d\in J\setminus
J^2$ such that $b=vdq^{n-1}$ and $q=dp^{\beta}$. Let $\gamma=max(\alpha,\beta)$. There exist $c'$ and $d'$ in
$J\setminus J^2$ such that $c=c'p^{\gamma -\alpha}$ and $d=d'p^{\gamma
  -\beta}$. It follows that $q=c'p^{\gamma}=d'p^{\gamma}$. Since $h$
is a monomorphism, we deduce that $c'=d'$.  Assume that
$m=n$. It follows that $b\in Ra$
if $\alpha\geq\beta$. Now assume that $n>m$. If $\beta\leq\alpha$ then
$b\in Ra$. If $\beta>\alpha$, let $q'\in J$ such that
$q=q'p^{\beta-\alpha}$. Then $b=vp^{\gamma -\alpha}c'q'q^{n-2}$ and
$b=vu^{-1}q'q^{n-m-1}a$. 

Now suppose that $J$ is nilpotent. Let $n$ be the least integer such
that $(JR_J)^n=0$. We easily deduce that $J^{n+1}=0$. We claim
that $J^n=
\ker h$. This is obvious if $n=1$. Assume that $n>1$. The inclusion
$J^n\subseteq \ker h$ is easy. Conversely let $0\not=r\in \ker h$. There
exists $c\in J\setminus J^2$ such that $Rr+Rq\subseteq Rc$. The
following equality $r=cd$ holds for some $d\in R$. Since $n>1$, $d\in
J$. If $d\in \ker h$
then there exists $s\in R\setminus J$ such that $sd=0$. But there exists
$c'\in J$ such that $c=sc'$. It follows that $r=dsc'=0$. Hence a
contradiction. Consequently $d\notin J^n$. Let $m$ be the greatest
integer such that $d\in J^m$. Then $m<n$,
$\displaystyle{\frac{d}{1}=\frac{q^m}{s'}}$ and
$\displaystyle{\frac{c}{1}=\frac{q}{t}}$ for some $s',t\in R\setminus
J$. It follows that
$\displaystyle{\frac{r}{1}=\frac{q^{m+1}}{s't}}$. We deduce that
$m=n-1$ and $r\in J^n$. Let $a$ and $b$ be two non-zero
elements of $J$. Then there exists $c\in J$ such that
$Ra+Rb+Rq=Rc$. It follows that $a=ca'$ and $b=cb'$ for some $a'$ and
$c'$ in $R$. Since $a\not=0$ and $b\not=0$, $a'$ and $b'$ are not in
$J^n$. From the first part of the proof we deduce
that $R/J^n$ is a valuation ring, so there
exist $r\in R$ and $s\in J^n$ such that $b'=ra'+s$. It follows that
$b=ra$.
\end{proof}

\begin{example}[Examples] \textnormal{If, in our Example~\ref{E:fg}, $D$ is a
  discrete valuation domain of Krull dimension$\geq 2$, then every artinian
  $R$-module is RD-injective. But this property is not satisfied if
  $D$ is a discrete rank one valuation domain. Ho\-wever, if $R$ is the
  ring defined in our Example~\ref{E:nfg}, then every artinian
  $R$-module is RD-injective.}

\textnormal{Let $R$ be the ring defined by the following pullback
  diagram of ring maps:
\[\begin{matrix}R&\rightarrow&T\\
\downarrow&{}&\downarrow\\
V&\rightarrow&K\end{matrix}\]
where $V$ is a discrete rank one valuation domain, $K$ its
  field of fractions and $T$ a local ring of residual class field
  $K$. Then $R$ is local and its maximal ideal $P$ is the inverse
  image of the maximal ideal of $V$. Clearly $R/P$ is RD-injective. If
  $J=R\cap Q$, where $Q$ is the maximal ideal of $T$, then $T=R_J$ and
  $V=R/J$. It is easy to prove that $R/J$ is RD-flat if and only if
  $T/Q$ is RD-flat over $T$.} 
\end{example}

We complete the proofs of Theorem~\ref{T:glo} and
Theorem~\ref{T:gloa}.

\begin{proof}[Proof of Theorem~\ref{T:glo}.]

$(1)\Leftrightarrow (7)$ by Theorem~\ref{T:main}.

Assume that $PR_P$ is not finitely generated over $R_P$
for every maximal ideal $P$. If the equivalent conditions (1) and (2) are
satisfied, then $P^2R_P=PR_P$ and $P=P^2$ for every maximal ideal
$P$. It follows that each artinian $R$-module is semi-simple of finite
length. Hence each artinian module is RD-flat and RD-injective. Let $M$ be a
noetherian $R$-module and $R'=R/\mathrm{ann}(M)$. Then the equality
$P^2R_P=PR_P$ implies that $R'_P=R/P$. It follows that $R'$ and $M$
are semi-simple.

Now we assume that there exists a maximal ideal $P$ such that $P\not= P^2$.

$(1)\Rightarrow (3)$ and $(2)\Rightarrow (4)$. Let $M$ be a module of finite
length. By
\cite[th\'eor\`eme p.368]{Bal76} there exists a finite family $F$ of maximal
ideals such that $M\cong \prod_{P\in F}M_P$. Consequently we may assume
that $M$ is indecomposable, $R$ is local and its
maximal ideal $P=pR$. There exists an integer $n$ such that $p^n$
annihilates $M$. Then it is easy to show that $M$ is RD-flat and
RD-injective.

$(5)\Rightarrow (6)$. Let $M$ be a noetherian module. By
\cite[lemme~1.3]{Cou82} the diagonal map $M\rightarrow\Pi_{P\in
  Max(R)}M_P$ is a pure monomorphism. For every maximal ideal $P$ we
put\\
\(F_P(M)=\mathrm{Hom}_{R_P}(\mathrm{Hom}_{R_P}(M_P,\mathrm{E}(R/P)),\mathrm{E}(R/P)).\) On
the other\\ hand, by \cite[Lemma XIII.2.3]{FuSa01}, $M_P$ is a pure submo\-dule of $F_P(M)$. Then $M$ is a pure submodule of $\Pi_{P\in
  \mathrm{Max}(R)}F_P(M)$. It is enough to prove that
  $F_P(M)$ is RD-injective. Consequently we may
assume that $R$ is local, $P$ is its maximal ideal and $E=\mathrm{E}(R/P)$. Let
$J=\cap_{n\in\mathbb{N}}Rp^n$, $E'=E[J]$ and $R'=R/J$. Thus $E'$ and
$\mathrm{Hom}_R(M,E)\cong \mathrm{Hom}_{R'}(M,E')$ are artinian. It
follows that
$\mathrm{Hom}_R(\mathrm{Hom}_R(M,E),E)$
is RD-injective.

Let $M$ be an indecomposable artinian module. By \cite[Th\'eor\`eme
p.368]{Bal76}, there exists a maximal ideal $P$ such that $M$ is an
$R_P$-module and every $R$-submodule of $M$ is also an
$R_P$-submodule. As in the proof of Theorem~\ref{T:artin} we prove
that $M$ is uniserial. 

If $M$ is a noetherian module then
$R'=R/\mathrm{ann}(M)$ is a noetherian RD-ring. Consequently $R'$ is a finite
product of Dedeking domains and artinian valuation rings. It follows
that $M$ is a direct sum of $2$-generated submodules.
\end{proof}

\bigskip
The following example shows that the converse of the last assertion of
Theorem~\ref{T:glo} is not true.
\begin{example}
\textnormal{Let $K$ be a field and for every $n\in\mathbb{N}$, let
  $R_n$ be the localization of the polynomial ring $K[X_n,Y_n]$ at the
  maximal ideal generated by $\{X_n,Y_n\}$ and $P_n$ the maximal
  ideal of $R_n$. For each $n\in\mathbb{N}$, let
  $\delta_n:R_n\rightarrow R_{n+1}$ be the ring homomorphism defined
  by $\delta_n(X_n)=X_{n+1}^2$ and $\delta_n(Y_n)=Y_{n+1}^2$. Let $R$
  be the direct limit of the system
  $(R_n,\delta_n)_{n\in\mathbb{N}}$. Then $R$ is local and its
  maximal ideal $P$ is the direct limit of
  $(P_n)_{n\in\mathbb{N}}$. It is obvious that $P^2=P$. Hence every
  artinian $R$-module is semi-simple of finite length. Let $r\in R$
  such that $X_0$ and $Y_0$ are in $Rr$. We may assume that there
  exist $n\in\mathbb{N}$, $s,t\in R_n$ such that $X_n^{2^n}=rs$ and
  $Y_n^{2^n}=rt$. It follows that $tX_n^{2^n}=sY_n^{2^n}$. Since $R_n$
  is a unique factorization domain there exists $u\in R_n$ such that
  $s=uX_n^{2^n}$ and $t=uY_n^{2^n}$. We deduce that
  $(1-ru)X_n^{2^n}=0$. Hence $r\notin P$ since $R$ is
  local. Consequently $R$ doesn't satisfy the equivalent conditions
  of Theorem~\ref{T:glo}.}
\end{example}

\begin{proof}[Proof of Theorem~\ref{T:gloa}.]

$(3)\Leftrightarrow (1)$. Let $M$ be an artinian module. By \cite[Th\'eor\`eme
p.368]{Bal76} we may assume that $R$ is local. Hence it is an immediate
consequence of Theorem~\ref{T:artin}.

$(1)\Rightarrow (2)$. By Proposition~\ref{P:flat} RD-flatness is a local
property, so we may assume that $R$ is local. Since $R/P$
is RD-flat, $R/J$ is noetherian and $E[J]=\mathrm{E}_{R/J}(R/P)$ is
artinian. Hence, if $M$ is a noetherian module then
\(\mathrm{Hom}_R(M,\mathrm{E}_R(R/P))=\mathrm{Hom}_{R/J}(M,E[J])\ \mathrm{is\
artinian.}\)
 We conclude by Proposition~\ref{P:flat}.

$(2)\Rightarrow (3)$. We may assume that $R$ is local. Since $R/P$ is
RD-flat, we may suppose that $P$ is principal. By
Corollary~\ref{C:cons}, $R/J$ is noetherian. By
Proposition~\ref{P:flat} it follows that $L(P)=E[J]$ is
RD-injective. 

$(3)\Leftrightarrow (4)$ by Theorem~\ref{T:artin}.
\end{proof}

\section{Complements to Section \ref{S:perf}}

This section is added to the version published in Communications in Algebra in 2006. Even if Theorem \ref{T:perf} is true, its proof is not complete. It is valid if we assume that $R$ is local. So, Corollaries \ref{C:local} and \ref{C:perf} are correctly proven. The aim of this added section is to give a complete proof of  the following theorem, where  \cite[Main Theorem]{Jon70} plays a crucial role.

\begin{theorem} \label{T:perf1}
Let $R$ be a ring. The following conditions are equivalent:
\begin{enumerate}
\item each RD-flat right $R$-module is RD-projective;
\item $R$ is right perfect.
\end{enumerate}
Moreover, each RD-projective right module is a direct sum of finitely presented cyclic indecomposable submodules if $R$ is right perfect.
\end{theorem}

The following lemma is needed to show this theorem. We shall use the same notations as in Section \ref{S:perf}.

\begin{lemma}\label{L:gen}
Let $R$ be a ring. For any element $A$ of $\mathcal{RD}$, $\mathrm{gen}\ A\leq 2$.
\end{lemma}
\begin{proof}
Let $X=R/A$. So, there exists $a\in R$ such that $R/aR=X\oplus Y$ where $Y$ is a submodule of $R/aR$. It is obvious that $X$ and $Y$ are cyclic. We consider the following commutative diagram with exact rows and columns, where the right vertical homomorphism is the identical map of $X$:

\[\begin{matrix}
 {} & {} & 0 & {} & 0 & {} & {} & {} & {}\\
{} & {} & \downarrow & {} & \downarrow & {} & {} & {} & {} \\
0 & \rightarrow & Z & \rightarrow & aR & \rightarrow & 0 & {} & {} \\
 {} & {} & \downarrow & {} & \downarrow & {} & \downarrow & {} & {} \\
0 & \rightarrow & A & \rightarrow & R & \rightarrow & X & \rightarrow & 0 \\
{} & {} & \downarrow & {} & \downarrow & {} & \downarrow & {} & {} \\
0 & \rightarrow & Y & \rightarrow & R/aR & \rightarrow & X & \rightarrow & 0 \\
{} & {} & \downarrow & {} & \downarrow & {} & \downarrow & {} & {} \\
{} & {} & 0 & {} & 0 & {} & 0 & {} & {} 
\end{matrix}\]

Since $Z\cong aR$, $Z$ is cyclic. From the the left vertical exact sequence we get that $\mathrm{gen}\ A\leq 2$.
\end{proof}

\begin{remark}
\textnormal{If $R$ is right perfect then $\mathcal{RD}$ verifies the ascending condition by \cite[Main Theorem]{Jon70}.}
\end{remark}

A ring $R$ is said to be \textbf{strongly $\pi$-regular} if, for each $r\in R$, there exist $s\in R$ and an integer $q\geq 1$ such that $r^q=r^{q+1}s$. By \cite[Theorem 3.16]{Fac98} each strongly $\pi$-regular $R$ satisfies the following condition: for each $r\in R$, there exist $s\in R$ and an integer $q\geq 1$ such that $r^q=sr^{q+1}$. Recall that a left  $R$-module $M$ is said to be {\bf Fitting} if for each endomorphism $f$ of $M$ there exists a positive integer $t$ such that $M=\ker\ f^t\oplus f^t(M)$.

\begin{proof}[Proof of Theorem \ref{T:perf1}]
By Theorem \ref{T:perf} and the beginning of its proof it remains to show 
$(2)\Rightarrow (1)$ and the last assertion.

We shall prove that $\mathcal{A}$ satifies the conditions (a) and (b) of the proof of Theorem \ref{T:perf}. For each positive integer $n$, the matrix ring $\mathrm{M}_n(R)$ is right perfect too. It follows that $\mathrm{M}_n(R)$ is strongly $\pi$-regular for each positive integer $n$. So, $\mathcal{A}$ satisfies the condition (a) by \cite[Theorem 3.22]{Fac98}. Moreover, by \cite[Lemmas 3.21 and 2.21]{Fac98}, $\mathrm{End}_R(X)$ is local for each finitely presented indecomposable object of $\mathcal{A}$. Since each right RD-projective module is a direct summand of a direct sum of finitely presented cyclic indecomposable modules, \cite[Theorems 2.8 and 2.9]{Fac98} imply the last assertion of our theorem.

Now, we consider the sequence \ref{E:seq} in the proof of Theorem \ref{T:perf}. For each integer $n>0$ let $A_n$ be the right ideal of $\mathcal{RD}$ for which $X_n\cong R/A_n$, and let $R_n=R$. For each integer $n>0$ there exist two homomorphisms $g_n:R_n\rightarrow R_{n+1}$ and $h_n:A_n\rightarrow A_{n+1}$ such that the following diagram commutes:

\[\begin{matrix}
A_n & \xrightarrow{h_n} & A_{n+1} \\
\downarrow & {} & \downarrow \\
R_n & \xrightarrow{g_n} & R_{n+1} \\
\downarrow & {} & \downarrow \\
X_n & \xrightarrow{f_n} & X_{n+1} 
\end{matrix}\]
For any integers $p>n>0$ we put $\alpha_{p,n}=f_p\circ\dots\circ f_n$, $\beta_{p,n}=g_p\circ\dots\circ g_n$ and $\gamma_{p,n}=h_p\circ\dots\circ h_n$. Let $(X,\ \alpha_n)$, $(R',\ \beta_n)$ and $(A,\ \gamma_n)$ be the direct limit of $\{X_n, \alpha_{p,n}\}$, $\{R_n, \beta_{p,n}\}$ and $\{A_n, \gamma_{p,n}\}$ respectively. Then the following sequence $0\rightarrow A\rightarrow R'\rightarrow X\rightarrow 0$ is exact. The family 
$(\alpha_n(X_n))_{n>0}$ is an ascending chain of cyclic submodules of $X$ whose union is $X$. By \cite[Main Theorem]{Jon70} there exists an integer $m>0$ such that 
$\alpha_n(X_n)=X$ for each integer $n\geq m$. In the same way we prove that $R'$ is cyclic. Since $R'$ is the direct limit of a system of free modules and $R$ is right perfect we deduce that $R'$ is projective. By using Lemma \ref{L:gen} and \cite[Main Theorem]{Jon70} we get that $\mathrm{gen}\ A\leq 2$. Hence, $X$ is finitely presented. By \cite[25.2(c)]{Wis91} the canonical homomorphism $\varinjlim_{n>0} \mathrm{Hom}_R(X,X_n)\rightarrow\mathrm{End}_R(X)$ is an isomorphism. So, there exist an integer $n>0$ and a homomorphism $\delta_n:X\rightarrow X_n$ such that $\alpha_n\circ\delta_n$ is the identical map of $X$. It follows that $X$ is isomorphic to a direct summand of $X_n$. If we assume that $X\ne 0$ then $\alpha_n$ is bijective because $X_n$ is indecomposable. In this case it is easy to see that 
$\alpha_p$ is an isomorphism for each integer $p\geq n$. This implies that $f_p$ is bijective for any $p\geq n$. We get a contradiction. Hence $X=0$, $\alpha_1=0$ and there exists an integer $N>1$ such that $\alpha_{N,1}=0$ because $X_1$ is cyclic. So, $f_N\circ\dots\circ f_2\circ f_1=0$ as required.
\end{proof}

\end{document}